\documentclass[10pt]{article}

\usepackage[longnamesfirst]{natbib}
\bibpunct{(}{)}{;}{a}{,}{,}
\setcitestyle{numbers,open={[},close={]}}

\usepackage[T1]{fontenc}
\usepackage[utf8]{inputenc}
\usepackage[english]{babel}
\usepackage{amsmath}
\usepackage{graphicx}
\usepackage{amssymb}
\usepackage{bm}
\usepackage{subdepth}
\usepackage{enumitem}
\usepackage{subcaption}
\usepackage{marvosym}

\usepackage[ruled,onelanguage]{algorithm2e}
\usepackage{algorithmic}

\usepackage[svgnames]{xcolor}
\usepackage{nameref}
\usepackage{hyperref}
\hypersetup{
  breaklinks,
  colorlinks,
  citecolor={blue},
  urlcolor={blue},
  linkcolor={blue}
}
\usepackage[nameinlink]{cleveref}

\makeatletter
\renewcommand*{\eqref}[1]{%
  \hyperref[{#1}]{\textup{\tagform@{\ref*{#1}}}}%
}
\makeatother

\newcommand{\dd}{\mathrm{d}}
\newcommand{\diff}[2]{\frac{\dd #1}{\dd #2}}
\newcommand{\T}{^\mathrm{T}}

\newcommand{\bv}{\ensuremath{\mathbf{b}}}
\newcommand{\pv}{\ensuremath{\mathbf{p}}}
\newcommand{\uv}{\ensuremath{\mathbf{u}}}
\newcommand{\vv}{\ensuremath{\mathbf{v}}}
\newcommand{\Fv}{\ensuremath{\mathbf{F}}}

\begin{document}
%
\title{An Optimal Control Problem\\
       for Elastic Registration and Force Estimation\\
       in Augmented Surgery}
\author{Guillaume Mestdagh\footnote{Inria, Strasbourg, France.
        \Letter\ \url{stephane.cotin@inria.fr}} 
        \and St\'{e}phane Cotin\footnotemark[1]}
\maketitle
\begin{abstract}
The nonrigid alignment between a pre-operative biomechanical model and an
intra-operative observation is a critical step to track the motion of a soft
organ in augmented surgery.
While many elastic registration procedures introduce artificial forces into the
direct physical model to drive the registration,
we propose in this paper a method to reconstruct the surface
loading that actually generated the observed deformation.
The registration problem is formulated as an optimal control problem where the
unknown is the surface force distribution that applies on the organ and the resulting
deformation is computed using an hyperelastic model.
Advantages of this approach include a greater control over the set of
admissible force distributions, in particular the opportunity to choose where
forces should apply, thus promoting physically-consistent displacement
fields.
The optimization problem is solved using a standard adjoint method.
We present registration results with experimental phantom data
showing that our procedure is competitive in terms of accuracy.
In an example of application, we estimate the forces applied by a surgery tool
on the organ.
Such an estimation is relevant in the context of robotic surgery systems, where
robotic arms usually do not allow force measurements, and providing force
feedback remains a challenge.
\end{abstract}

\paragraph{Keywords}
Augmented surgery, Optimal control, Biomechanical simulation

\section{Introduction}
\label[section]{intro}
In the context of minimally-invasive surgery, abdominal organs are constantly
subject to large deformations.
These deformations, in combination with the limited visual feedback available
to medical staff, make an intervention such as tumor ablation a complex task
for the surgeon.
Augmented reality systems have been designed to provide a three-dimensional
view of an organ, which shows the current position of internal structures.
This virtual view is superimposed onto intra-operative images displayed in the
operating room.

Without loss of generality, we consider the case of liver laparoscopic surgery.
In \cite{haouchine2013}, the authors describe a full pipeline to produce
augmented images.
In their work, available data are intra-operative images provided by
a laparoscopic stereo camera, and a biomechanical model of the organ and its
internal structures, computed from pre-operative CT scans.
During the procedure, a point cloud representing the current location of the
liver surface is first extracted from laparoscopic images, and then the
pre-operative model is aligned with the intra-operative point cloud in a
non-rigid way.

The elastic registration procedure used to perform the alignment between
pre- and intra-operative data is a widely studied subject.
In particular, much attention has been given to choosing an accurate
direct mechanical model for the liver and its surroundings.
The liver parenchyma is usually described using a hyperelastic constitutive
law (see \cite{marchesseau2017} for an extensive review).
Used hyperelastic models include Saint-Venant-Kirchhoff \cite{delingette2004},
neo-Hookean \cite{miller2007} or Ogden models~\cite{nikolaev2020}.
Due to its reduced computational cost, the linear co-rotational model~\cite{nesme2005}
is a popular solution when it comes to matching a real-time
performance constraint \cite{suwelack2014,plantefeve2016,peterlik2018}.
Additional stiffness due to the presence of blood vessels across the parenchyma
is sometimes also taken into account \cite{haouchine2015}.
A less discussed aspect of the modeling is the interaction between the liver
and its surroundings, which results in boundary conditions applied on the
liver surface.
Proposed approaches often involve Dirichlet boundary conditions where the main
blood vessels enter the liver \cite{peterlik2018} or springs 
to represent ligaments holding the organ \cite{nikolaev2020}.

Another key ingredient to obtain accurate reconstructions is the choice of
a registration procedure.
In many physics-based registration methods in the literature, fictive forces or
energies are added into the direct problem to drive the registration.
Approaches based on the Iterative Closest Point algorithm introduce attractive
forces beetween the liver and the observed point cloud and let the system evolve
as a time-dependent process toward an equilibrium
\cite{haouchine2013,plantefeve2016}.
In~\cite{suwelack2014}, the authors model those attractive forces using an
electrostatic potential.
In \cite{peterlik2018}, sliding constraints are
used in the direct problem to enforce correspondence between the deformed
liver model and the observed data.
Such constraints are enforced by Lagrange multipliers that are also fictive
forces.
As the intra-operative observation is not a real protagonist of the physical
model, those fictive forces do not reflect the true causes of displacements.
This results in a poorly physically-consistent displacement field, regardless
of the direct model.
These methods cannot guarantee accurate registrations.

In this paper, we present a method to reconstruct a
surface force distribution that explains the observed deformation, using the 
optimal control formalism.
The problem formulation includes the choice of an elastic model
(\Cref{sec:21}), a set of admissible forces and an objective function to
minimize (\Cref{sec:22}).
In particular, the set of admissible forces is specified by the direct model
to ensure physically-consistent deformations.
The adjoint method used to perform the optimization is described in
\Cref{sec:23}.

Some existing works in the literature are concerned with reconstructing the
physical causes of displacement.
In \cite{ozgur2018}, the effects of gravity and pneumoperitoneum pressure are
taken into account to perform an initial intraoperative registration.
In \cite{rucker2014}, the authors control the imposed displacement on parts of
the liver boundary that are subject to contact forces, while a free boundary
condition is applied onto other parts.
A numerical method similar to our adjoint method is used in 
\cite{heiselman2020} to register a liver model onto a point cloud coupled with
ultrasound data.
While their method is very specific and tailored for linear elasticity, we
present a more generic approach which is compatible with nonlinear models.
Our method is relevant when a measurement of forces is needed, in particular
when the surgeon interacts with the organ through a surgical robotic system
without force sensor.
We give an example of such a force estimation in \Cref{sec:32}.

\section{Methods}
\label[section]{sec:2}
In this section, we give details about the optimal control problem formulation.
We first specify some notation around the direct model, then we introduce the
optimization problem.
In the last part, we describe the adjoint method, which is used to compute
derivatives of the objective function with respect to the control.

\subsection{Hyperelastic model and observed data}
\label[section]{sec:21}
The liver parenchyma in its pre-operative configuration is 
represented by a meshed domain $\Omega$, filled with an elastic material.
When a displacement field $\uv$ is applied to $\Omega$, the deformed mesh is
denoted by $\Omega_\uv$ and its boundary is denoted by $\partial \Omega_\uv$.
Note that the system state is fully known through the displacement of mesh
nodes, stored in $\uv$.
The liver is embedded with a hyperelastic model.
When a surface force distribution $\bv$ is applied to the liver boundary, the 
resulting displacement $\uv_\bv$ is the unique solution of the static
equilibrium equation
\begin{equation}\label{eq:hyperelastic}
    \Fv(\uv_\bv) = \bv,
\end{equation}
where $\Fv$ is the residual from the hyperelastic model.
Note that $\Fv$ is very generic and may also account for blood
vessels rigidity, gravity or other elements in a more sophisticated
direct model.

The observed intra-operative surface is represented by a point cloud
$\Gamma = \{y_1, \dots, y_m\}$.
We also define the orthogonal projection onto $\partial\Omega_\uv$, also called
the closest point operator,
\begin{equation*}
    p_\uv(y) = \arg\min_{x\in \partial \Omega_\uv} \|y - x\|.
\end{equation*}
Here, we consider that the orthogonal projection operator always returns a
unique point, as points with multiple projections onto
$\partial\Omega_\uv$ represent a negligible subset of $\mathbb{R}^3$.

\subsection{Optimization problem}
\label[section]{sec:22}
We perform the registration by computing a control $\bv$ so that the
resulting displacement $\uv_\bv$ is in adequation with observed data.
The optimization problem reads
\begin{equation}\label{eq:opt-problem}
    \min_{\bv\in B}\quad \Phi(\bv)
    \qquad\text{where}\qquad
    \Phi(\bv) = J(\uv_\bv) + R(\bv),
\end{equation}
where $J$ measures the discrepancy between a prediction
$\uv_\bv$ and observed data, $R$ is an optional regularization term, and $B$ is
the set of admissible controls.

The functional $J$ enforces adequation with observed data and contains 
information such as corresponding landmarks between pre-- and
intra-operative surfaces.
Many functionals of this kind exist in the literature, based on
surface correspondence tools
\cite[see][and references therein]{sahillioglu2020}.
In this paper, we use a simple least-squares term involving the orthogonal
projection onto $\partial \Omega_\uv$, which reads
\begin{equation}\label{eq:functional}
    J(\uv) = \tfrac{1}{2m} \sum_{j=1}^m \|p_\uv(y_j) - y_j\|^2.
\end{equation}
The functional $J(\uv)$ evaluates to zero whenever every point $y\in\Gamma$
is matched by the deformed surface $\partial\Omega_\uv$.
If $\vv$ is a perturbation of the current displacement field $\uv$, the
gradient of $\eqref{eq:functional}$ with respect to the displacement is defined
by
\begin{equation}
    \left\langle \nabla J(\uv), \vv \right\rangle
    = \tfrac{1}{m} \sum_{j=1}^m \left(r_j - y_j\right)
    \cdot \vv(r_j),
\end{equation}
where $r_j = p_\uv(y_j)$ and $\vv(r_j)$ is the displacement of the point
of $\partial\Omega_\uv$ currently at $r_j$ under the perturbation $\vv$.

The set of admissible controls $B$ contains a priori information about the
surface force distribution $\bv$, 
including the parts of the liver boundary where $\bv$ is nonzero and the
maximal intensity it is allowed to take.
The selection of zones where surface forces apply is critical to obtain
physically plausible registrations. 
In comparison, a spring-based approach would result in forces concentrated in
zones where the spring are fixed, which might disagree with the direct
model.

\subsection{Adjoint method}
\label[section]{sec:23}
We solve problem \eqref{eq:opt-problem} using an adjoint method.
In such a method, the only optimization variable is $\bv$.
Descent directions for the objective function $\Phi$ are computed in the space
of controls, which requires to compute $\nabla\Phi(\bv)$.
To differentiate $J(\uv_\bv)$ with respect to $\bv$, we use an adjoint
state $\pv_\bv$ defined as the solution of the adjoint system
\begin{equation}\label{eq:adjoint-equation}
    \nabla\Fv(\uv_\bv)\T \pv_\bv = \nabla J(\uv_\bv).
\end{equation}
A standard calculation \cite[see for instance][]{allaire2007} shows that
\begin{equation*}
    \pv_\bv = \diff{}{\bv}\left[ J(\uv_\bv)\right]
    \qquad\text{and}\qquad
    \nabla \Phi(\bv) = \pv + \nabla R(\bv).
\end{equation*}

The procedure to compute the objective gradient is summarized in
\Cref{alg:adjoint}.
For a given $\bv$, evaluating $\Phi(\bv)$ and $\nabla\Phi(\bv)$ requires to
solve the direct problem, and then to assemble and solve the adjoint
problem, which is linear.
In other words, the additional cost compared to a direct simulation is that
of solving one linear system.
The resulting gradient is then fed to a standard gradient-based optimization
algorithm to solve \eqref{eq:opt-problem} iteratively.

\begin{algorithm}
    \SetAlgoLined
    \KwData{Current iterate $\bv$}
    Compute the displacement $\uv_\bv$ by solving \eqref{eq:hyperelastic}\\
    Evaluate $\nabla J(\uv_\bv)$ and $\nabla \Fv(\uv_\bv)$\\
    Compute the adjoint state $\pv_\bv$ by solving \eqref{eq:adjoint-equation}\\
    \KwResult{$\nabla \Phi(\bv) = \pv_\bv + \nabla R(\bv)$}
    \label[algorithm]{alg:adjoint}
    \caption{Computation of objective gradient using an adjoint method.}
\end{algorithm}

\section{Results}
\label[section]{sec:3}
Our method is implemented in Python (Numpy and Scipy), and we use a 
limited-memory BFGS algorithm \cite{byrd1995} as an optimization procedure.
Our numerical tests run on an Intel Core i7-8700 CPU at 3.20~GHz with
16~GB RAM.

\subsection{Sparse Data Challenge dataset}
\label[section]{sec:31}
To evaluate our approach in terms of displacement accuracy,
we use the Sparse Data Challenge\footnote{See details and results at 
\url{sparsedatachallenge.org}.} dataset.
It consists of one tetrahedral mesh representing a liver phantom in its initial
configuration and 112 point clouds acquired from deformed configurations of
the same phantom \cite{collins2017,brewer2019}.
Once the registration is done, the final position of the mesh nodes for each
point cloud is uploaded on the challenge website.
Then the target registration error (TRE) is computed, using 159 targets
whose positions are unknown to us.

\begin{figure}
    \centering
    \includegraphics[width=.39\textwidth]{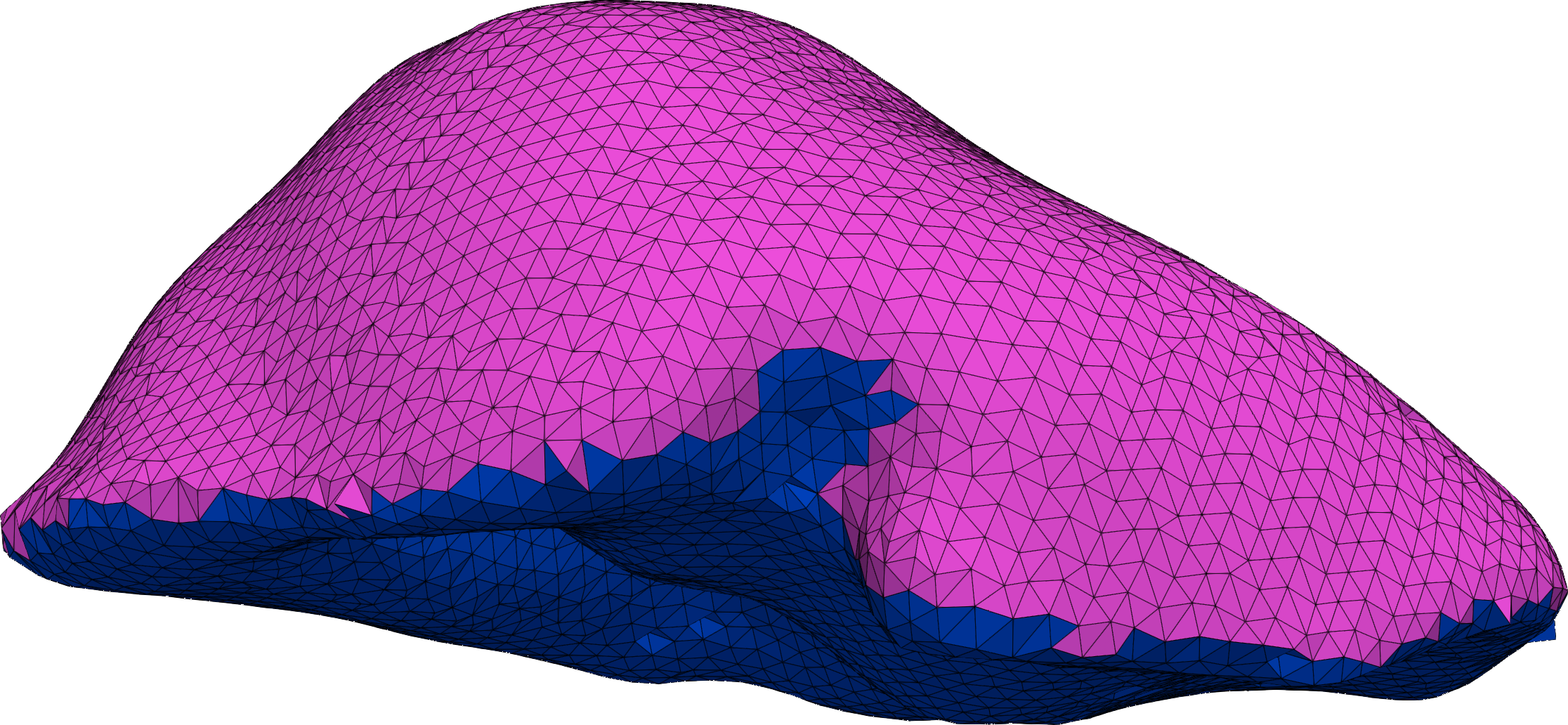}
    \includegraphics[height=.24\textwidth]{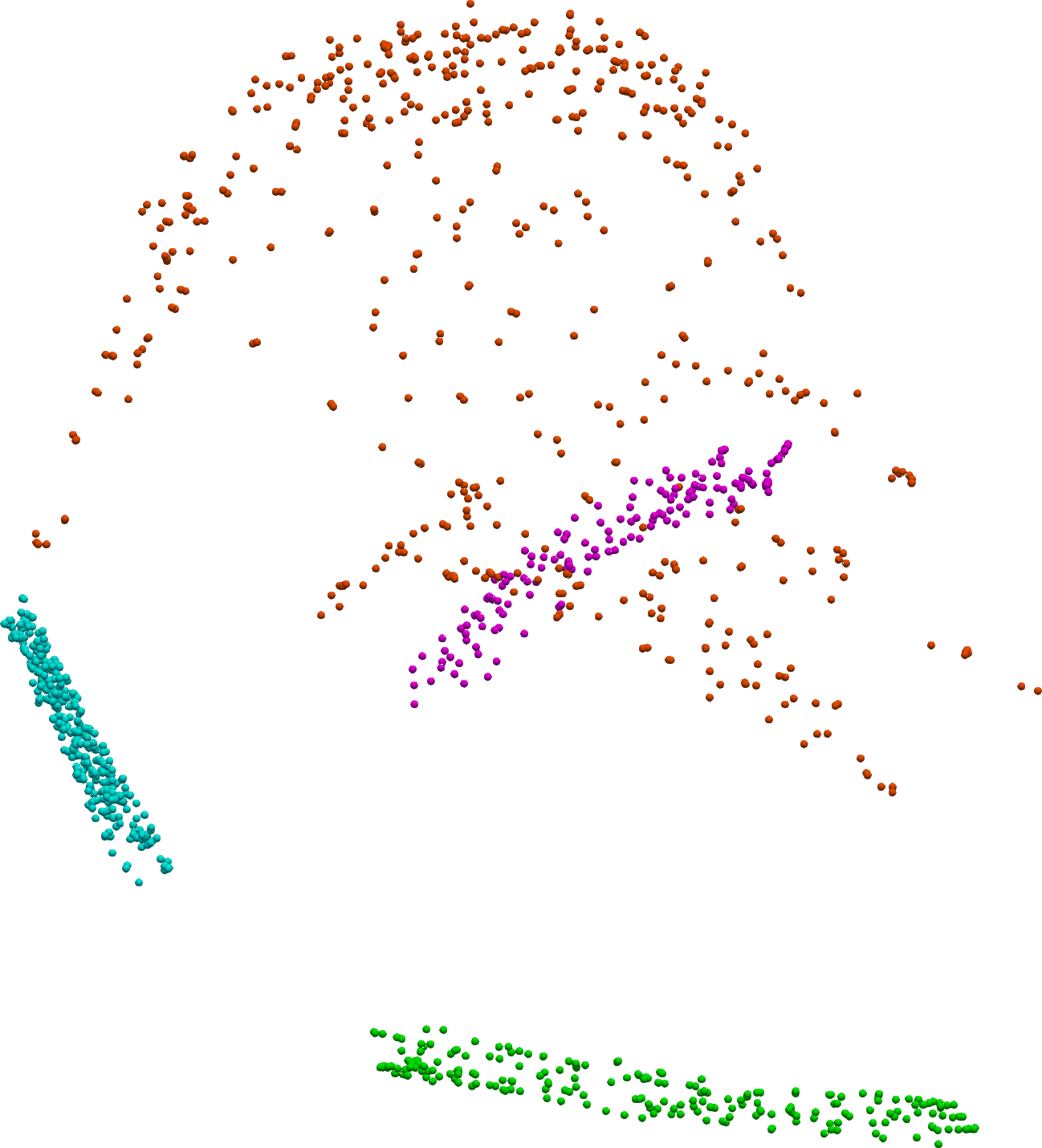}
    \includegraphics[height=.24\textwidth]{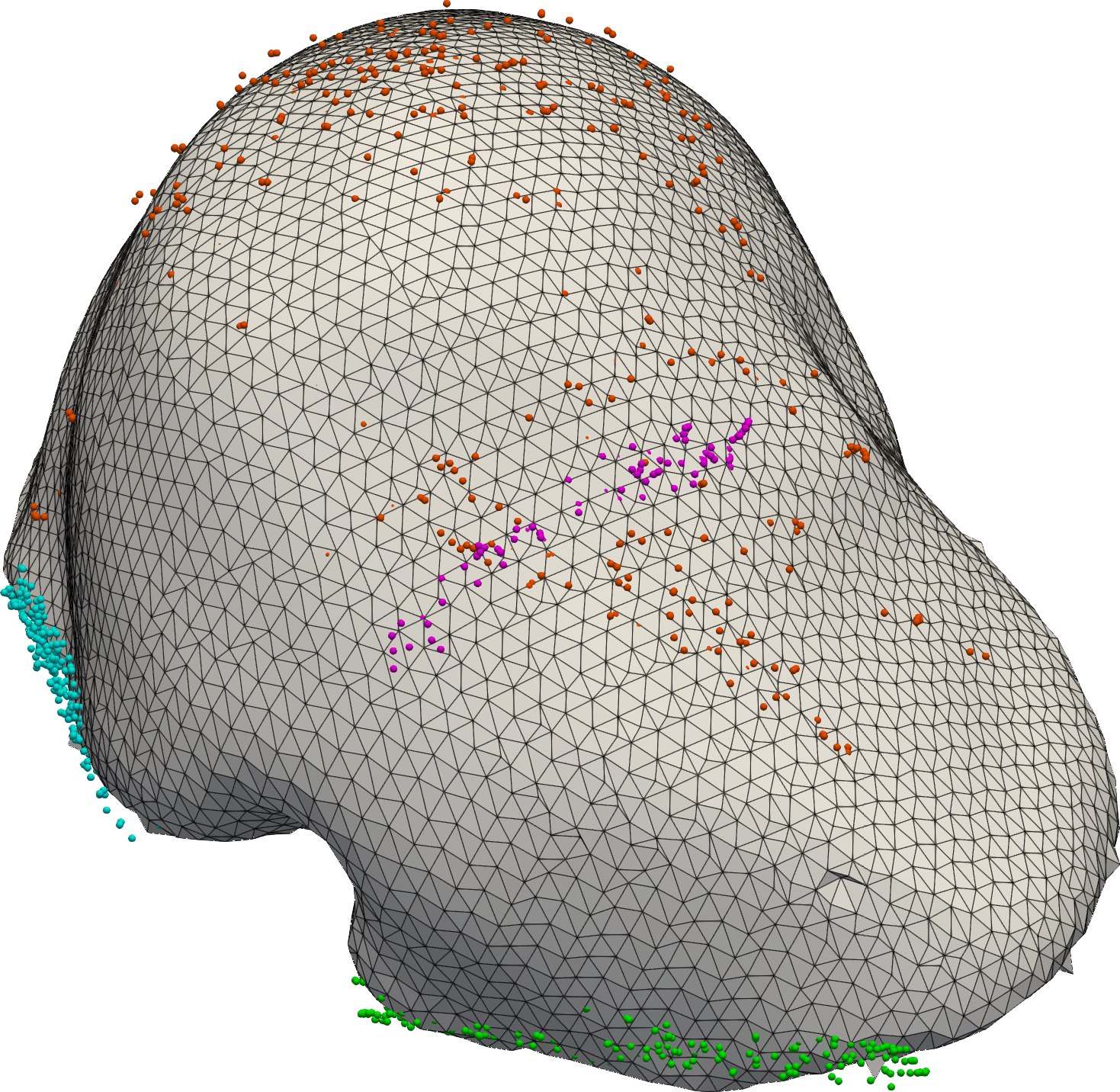}
    \caption{Left: Initial liver mesh with matching surface on top and loaded
             surface in the bottom. Right: Point cloud and liver mesh after
             rigid registration}
    \label[figure]{fig:sdc-dataset}
\end{figure}

Before we begin the elastic registration process, we use the standard
Iterative Closest Point method \cite{besl1992} to perform a rigid alignment.
Then we set a fixed boundary condition on a small zone of the posterior face
to enforce the uniqueness of the solution to the direct elastic problem.
As the liver mesh represents a phantom, no information about blood vessels is
available and for this reason we just choose six adjacent triangles close to
the center of the posterior face to apply the fixed displacement constraint.
As specified by the challenge authors, the forces causing the deformation are 
contact forces applied onto the posterior face of the phantom, while the point
cloud was acquired by observing the anterior surface.
As a consequence, we label the anterior surface as the "matching surface",
which is to be matched with the point cloud, while the posterior surface is
labeled as "loaded surface", where the force distribution $\bv$ is allowed
to take nonzero values.
\Cref{fig:sdc-dataset} shows the initial mesh with the matching and loaded
surfaces in different colors, a point cloud and the liver mesh after rigid
registration.
In a clinical context, the surface in the field of view of the
camera may be labeled as the matching surface while forces are allowed on the
remaining hidden surface.

We solve problem \eqref{eq:opt-problem} (with $R(\bv) = 0$) using the adjoint
method presented in \Cref{sec:23}.
As only small deformations are involved in this dataset, we use a linear
elastic model for the liver phantom, with $E = 1$ and $\nu = 0.4$.
The procedure is stopped after 200 iterations.

\begin{figure}
    \centering
    \includegraphics[width=\textwidth]{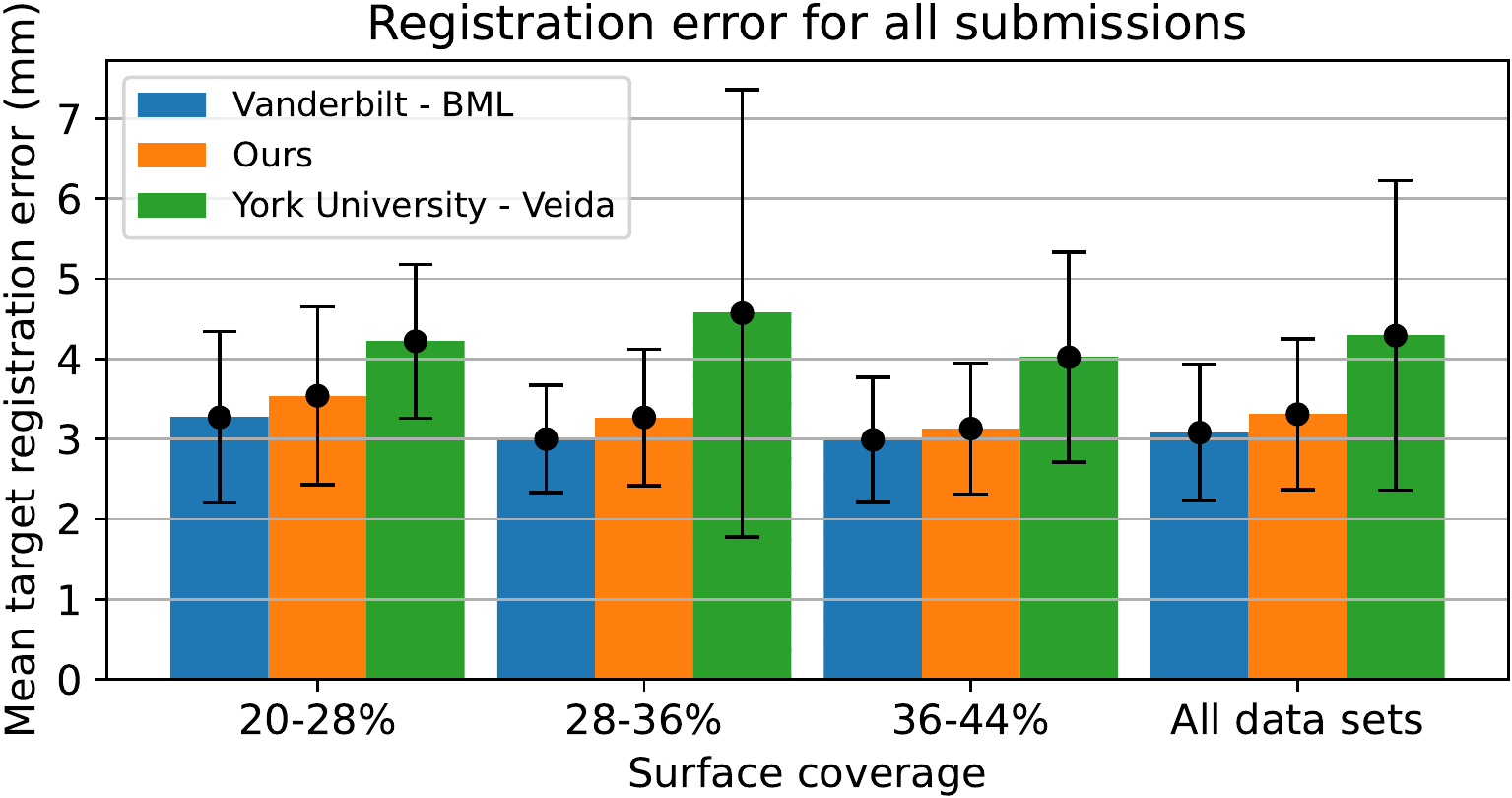}
    \caption{Comparison of our TRE results with other participants to the
             challenge.}
    \label[figure]{fig:sdc-bars}
\end{figure}

In \Cref{tab:sdc-results}, we reported the target registration error
statistics for all datasets returned by the challenge website after
submission of the results.
Even for point clouds with a low surface coverage, the average target
registration error stays below the 5~mm error that is usually required for
clinical applications.
In \Cref{fig:sdc-bars}, our results are compared with other submissions to the
challenge.
We obtained the second best result displayed on the challenge website (the
first one is that of the challenge organizers), which shows that our
approach is competitive for registration applications.
No information about methods used by other teams is displayed on the
challenge website.

\begin{table}
    \centering
    \begin{tabular}{|c|c|c|c|}
        \hline
        Surface Coverage & Average & Standard deviation & Median\\
        \hline
        20-28~\%          & 3.54    & 1.11               & 3.47\\
        28-36~\%          & 3.27    & 0.85               & 3.19\\
        36-44~\%          & 3.13    & 0.82               & 3.13\\
        All data sets    & 3.31    & 0.94               & 3.19\\
        \hline
    \end{tabular}
    \caption{Target registration error statistics (in millimeters) for all
             datasets, as returned by the website after submission.}
    \label[table]{tab:sdc-results}
\end{table}

\subsection{Force estimation in robotic surgery}
\label[section]{sec:32}
Estimating the force applied by a robotic arm onto an organ is necessary to
avoid causing damage to living tissues.
As certain standard surgical robots are not equipped with force sensors,
many indirect methods based on image processing have been proposed
\cite{nazari2021}.
Here we estimate a force in a context similar to
\cite{haouchine2018}, where the intra-operative point cloud is estimated from a 
laparoscopic camera.

We generate 5 synthetic test cases using a liver mesh of 3,046 vertices and a 
linear elastic model ($E = 20,000\ \mathrm{Pa}$, $\nu = 0.45$).
Dirichlet conditions are applied at the hepatic vein entry and in the falciform
region.
A test case is a sequence of traction forces applied onto two
adjacent triangles on the anterior surface of the liver to mimick the action
of a robotic tool manipulating the liver.
For each traction force, the resulting displacement is computed and a part of 
the deformed boundary is sampled to create a point cloud of 500 points.
Each sequence consists of 50 forces, with a displacement of about 
1~mm between two successive forces.

\begin{figure}
    \centering
    \includegraphics[width=.3\textwidth]{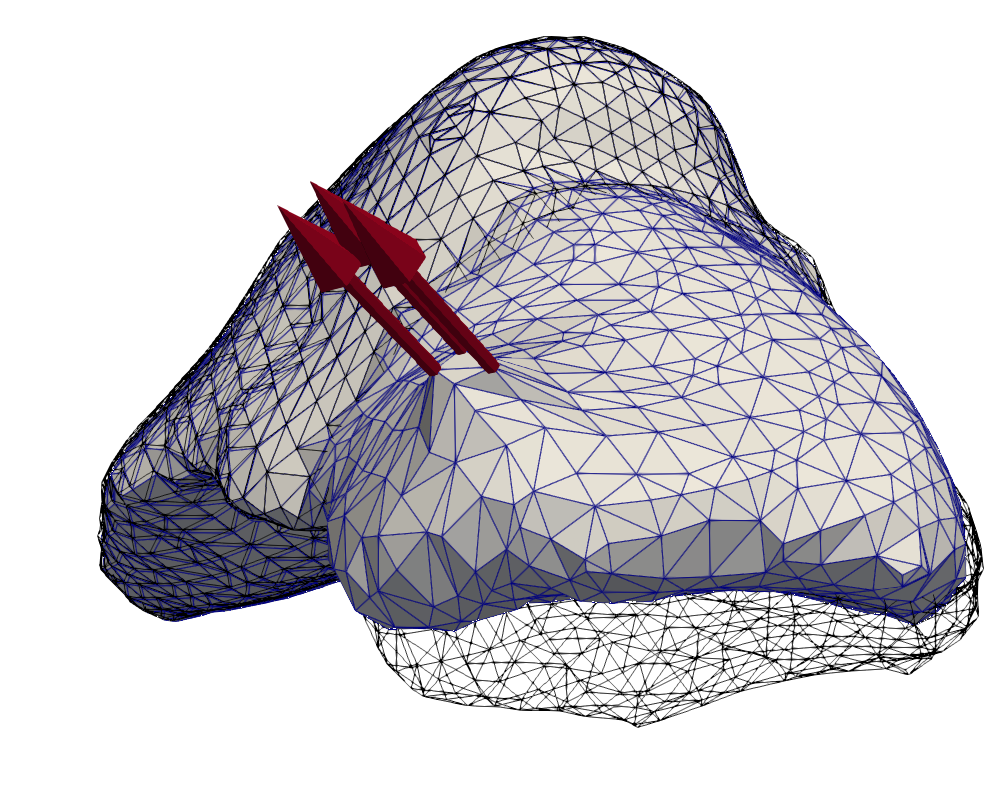}
    \includegraphics[width=.3\textwidth]{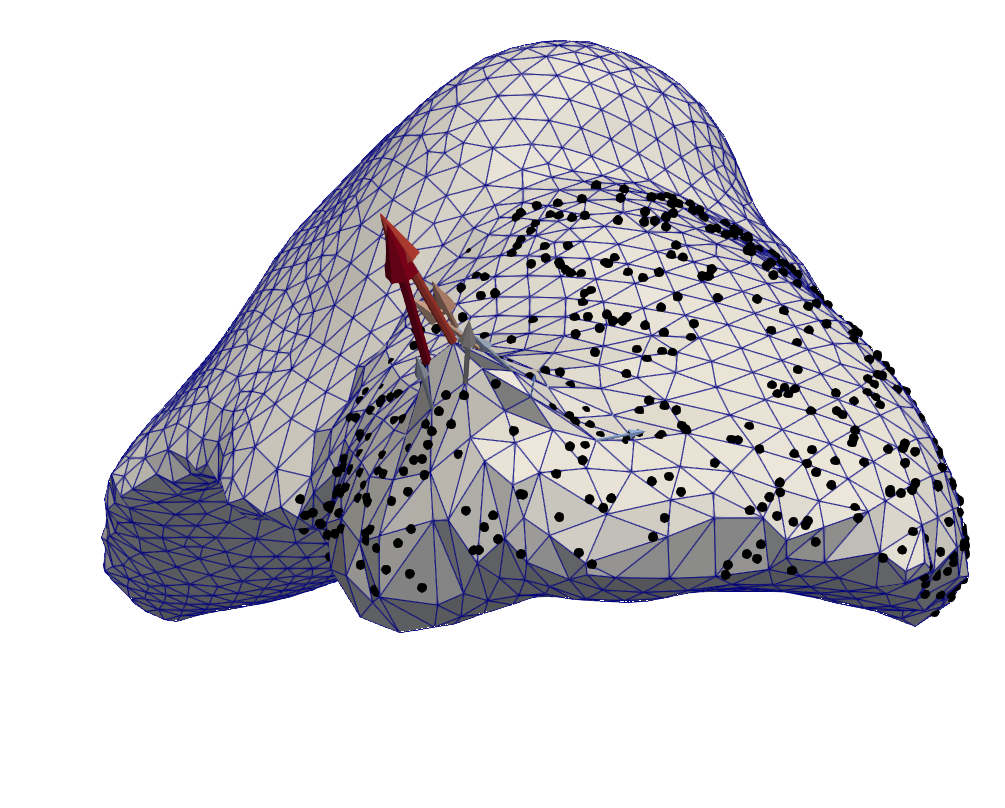}
    \includegraphics[width=.3\textwidth]{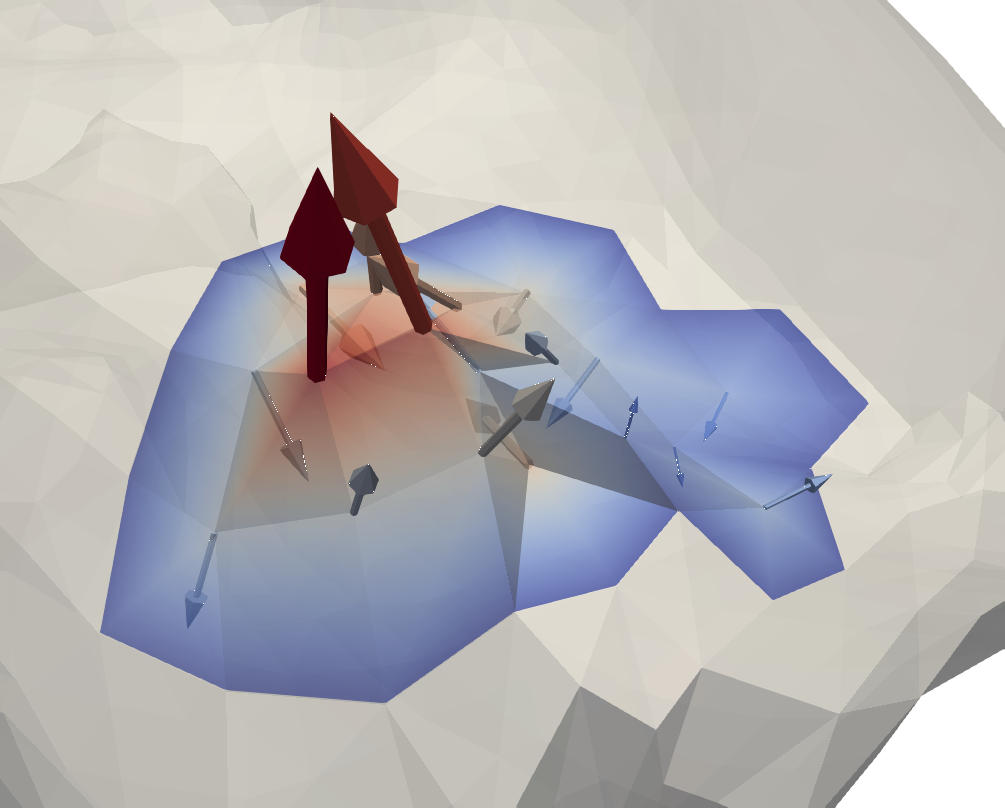}
    \caption{Synthetic deformation generated by a local force (left), 
             reconstructed deformation using the point cloud (center) and 
             zoom on the reconstructed force field (right).
             Nodal forces are summed to produce a resultant estimation.}
    \label[figure]{fig:forces}
\end{figure}

To avoid the inverse crime, we use a different mesh (3,829 vertices) for the
registration, with the same linear elastic model as above.
We allow the force distribution $\bv$ to be nonzero
only in a small zone (about 50 vertices) surrounding the triangles
concerned by the traction force.
In a clinical context, this approximate contact zone may be determined by
segmenting the instrument tip on laparoscopic images.
For a given traction force $f_\mathrm{true}$, we solve the optimization problem
with a relative tolerance of $5\cdot 10^{-4}$ on
the objective gradient norm and we compute the force estimation
$f_\mathrm{est}\in\mathbb{R}^3$ as the resultant of all the nodal
forces of the reconstructed distribution~$\bv$.
Then we compute the relative error 
$\|f_\mathrm{est} - f_\mathrm{true}\| / \|f_\mathrm{true}\|$.
\Cref{fig:forces} shows the original synthetic deformation and the
reconstructed deformation and surface force distribution.

For each sequence, the observed point clouds are successively fed to the
procedure to update the force estimation by solving the optimization
problem.
The optimization algorithm is initialized with the last reconstructed distribution, 
so that only a few iterations are required for the update.
\Cref{tab:force-results} shows the average error obtained for each sequence,
together with the average execution time of the updates and the number of
evaluations of the objective function.

The noise in the reconstructed distribution results in an intensity
smaller than the reference, which is the main cause of an overall 
error of 10~\% in average.
Note that the estimated force is proportional to the Young modulus of the model,
which can be measured using elastography.
According to \cite{oudry2014}, an error of 20~\% is typical for clinical
elastographic measurements.
In this context, the elastographic estimation is responsible for a larger part of the 
total error on the force estimation than our method.

\begin{table}
    \centering
    \begin{tabular}{|c|c|c|c|c|}
        \hline
        Sequence & N. eval. & Update Time & Relative error\\
        \hline
        Case 1 & 9.2 & 1.42~s &  8.9~\%  \\
        Case 2 & 5.6 & 0.85~s & 16.2~\%  \\
        Case 3 & 5.5 & 0.84~s &  5.7~\%  \\
        Case 4 & 5.4 & 0.82~s &  4.4~\%  \\
        Case 5 & 6.2 & 0.96~s & 15.0~\% \\
        \hline
    \end{tabular}
    \caption{Average number of objective evaluations, execution time per update 
             and relative error for each sequence.}
    \label[table]{tab:force-results}
\end{table}

\section{Conclusion}
\label[section]{conclu}
We presented a formulation of the liver registration problem using the
optimal control formalism, and used it with success in two different
application cases.
We showed that we can not only reconstruct an accurate displacement field,
but also, in certain cases, give a meaning to the optimal force distribution
returned by the procedure.
By tuning the formulation parameters (namely the set of admissible controls), 
we easily added new hypotheses into the direct model without changing the code
of our procedure.
These numerical results highlight the relevance and the flexibility of an 
optimal control approach in augmented surgery.
Due to its mathematical foundations, the optimal control formalism is
probably an important step toward provable accuracy for registration methods.
Current limitations in our work include the lack of results with
hyperelastic models due to their high computational cost. 
In the next steps of our work, we intend to reduce computation times by
initializing the optimization process with the output of a neural network.


\end{document}